\documentclass[reqno,12pt,a4paper]{amsart}

\usepackage{mathrsfs}
\usepackage{amssymb}
\usepackage{amsfonts}
\usepackage{latexsym}
\usepackage{amsthm}
\usepackage{amsmath}
\usepackage[arrow, matrix, curve]{xy}
\usepackage{graphicx}
\usepackage{tikz}
\usetikzlibrary{matrix,arrows}
\input xy

\begin{document}


\renewcommand{\theequation}{\arabic{section}.\arabic{equation}}
\theoremstyle{plain}
\newtheorem{theorem}{\bf Theorem}[section]
\newtheorem{lemma}[theorem]{\bf Lemma}
\newtheorem{corollary}[theorem]{\bf Corollary}
\newtheorem{proposition}[theorem]{\bf Proposition}
\newtheorem{definition}[theorem]{\bf Definition}
\newtheorem{remark}[theorem]{\it Remark}

\def\a{\alpha}  \def\cA{{\mathcal A}}     \def\bA{{\bf A}}  \def\mA{{\mathscr A}}
\def\b{\beta}   \def\cB{{\mathcal B}}     \def\bB{{\bf B}}  \def\mB{{\mathscr B}}
\def\g{\gamma}  \def\cC{{\mathcal C}}     \def\bC{{\bf C}}  \def\mC{{\mathscr C}}
\def\G{\Gamma}  \def\cD{{\mathcal D}}     \def\bD{{\bf D}}  \def\mD{{\mathscr D}}
\def\d{\delta}  \def\cE{{\mathcal E}}     \def\bE{{\bf E}}  \def\mE{{\mathscr E}}
\def\D{\Delta}  \def\cF{{\mathcal F}}     \def\bF{{\bf F}}  \def\mF{{\mathscr F}}
\def\c{\chi}    \def\cG{{\mathcal G}}     \def\bG{{\bf G}}  \def\mG{{\mathscr G}}
\def\z{\zeta}   \def\cH{{\mathcal H}}     \def\bH{{\bf H}}  \def\mH{{\mathscr H}}
\def\e{\eta}    \def\cI{{\mathcal I}}     \def\bI{{\bf I}}  \def\mI{{\mathscr I}}
\def\p{\psi}    \def\cJ{{\mathcal J}}     \def\bJ{{\bf J}}  \def\mJ{{\mathscr J}}
\def\vT{\Theta} \def\cK{{\mathcal K}}     \def\bK{{\bf K}}  \def\mK{{\mathscr K}}
\def\k{\kappa}  \def\cL{{\mathcal L}}     \def\bL{{\bf L}}  \def\mL{{\mathscr L}}
\def\l{\lambda} \def\cM{{\mathcal M}}     \def\bM{{\bf M}}  \def\mM{{\mathscr M}}
\def\L{\Lambda} \def\cN{{\mathcal N}}     \def\bN{{\bf N}}  \def\mN{{\mathscr N}}
\def\m{\mu}     \def\cO{{\mathcal O}}     \def\bO{{\bf O}}  \def\mO{{\mathscr O}}
\def\n{\nu}     \def\cP{{\mathcal P}}     \def\bP{{\bf P}}  \def\mP{{\mathscr P}}
\def\r{\rho}    \def\cQ{{\mathcal Q}}     \def\bQ{{\bf Q}}  \def\mQ{{\mathscr Q}}

 \def\cR{{\mathcal R}}     \def\bR{{\bf R}}  \def\mR{{\mathscr R}}
                \def\cS{{\mathcal S}}     \def\bS{{\bf S}}  \def\mS{{\mathscr S}}
\def\t{\tau}    \def\cT{{\mathcal T}}     \def\bT{{\bf T}}  \def\mT{{\mathscr T}}
\def\f{\phi}    \def\cU{{\mathcal U}}     \def\bU{{\bf U}}  \def\mU{{\mathscr U}}
\def\F{\Phi}    \def\cV{{\mathcal V}}     \def\bV{{\bf V}}  \def\mV{{\mathscr V}}
\def\P{\Psi}    \def\cW{{\mathcal W}}     \def\bW{{\bf W}}  \def\mW{{\mathscr W}}
\def\o{\omega}  \def\cX{{\mathcal X}}     \def\bX{{\bf X}}  \def\mX{{\mathscr X}}
\def\x{\xi}     \def\cY{{\mathcal Y}}     \def\bY{{\bf Y}}  \def\mY{{\mathscr Y}}
\def\X{\Xi}     \def\cZ{{\mathcal Z}} 
\def\cm{{\mathcal m}}
 \def\bZ{{\bf Z}}  \def\mZ{{\mathscr Z}}
\def\O{\Omega}

\newcommand{\gA}{\mathfrak{A}}
\newcommand{\gB}{\mathfrak{B}}
\newcommand{\gC}{\mathfrak{C}}
\newcommand{\gD}{\mathfrak{D}}
\newcommand{\gE}{\mathfrak{E}}
\newcommand{\gF}{\mathfrak{F}}
\newcommand{\gG}{\mathfrak{G}}	
\newcommand{\gH}{\mathfrak{H}}
\newcommand{\gI}{\mathfrak{I}}
\newcommand{\gJ}{\mathfrak{J}}
\newcommand{\gK}{\mathfrak{K}}
\newcommand{\gL}{\mathfrak{L}}
\newcommand{\gM}{\mathfrak{M}}
\newcommand{\gN}{\mathfrak{N}}
\newcommand{\gO}{\mathfrak{O}}
\newcommand{\gP}{\mathfrak{P}}
\newcommand{\gQ}{\mathfrak{Q}}
\newcommand{\gR}{\mathfrak{R}}
\newcommand{\gS}{\mathfrak{S}}
\newcommand{\gT}{\mathfrak{T}}
\newcommand{\gU}{\mathfrak{U}}
\newcommand{\gV}{\mathfrak{V}}
\newcommand{\gW}{\mathfrak{W}}
\newcommand{\gX}{\mathfrak{X}}
\newcommand{\gY}{\mathfrak{Y}}
\newcommand{\gZ}{\mathfrak{Z}}

\def\ve{\varepsilon} \def\vt{\vartheta} \def\vp{\varphi}  \def\vk{\varkappa}

\def\Z{{\mathbb Z}} \def\R{{\mathbb R}} \def\C{{\mathbb C}}  \def\K{{\mathbb K}}
\def\T{{\mathbb T}} \def\N{{\mathbb N}} \def\dD{{\mathbb D}} 
\def\B{{\mathbb B}}


\def\la{\leftarrow}              \def\ra{\rightarrow}     \def\Ra{\Rightarrow}
\def\ua{\uparrow}                \def\da{\downarrow}
\def\lra{\leftrightarrow}        \def\Lra{\Leftrightarrow}
\newcommand{\abs}[1]{\lvert#1\rvert}
\newcommand{\br}[1]{\left(#1\right)}

\def\lan{\langle} \def\ran{\rangle}


\def\lt{\biggl}                  \def\rt{\biggr}
\def\ol{\overline}               \def\wt{\widetilde}
\def\no{\noindent}


\let\ge\geqslant                 \let\le\leqslant
\def\lan{\langle}                \def\ran{\rangle}
\def\/{\over}                    \def\iy{\infty}
\def\sm{\setminus}               \def\es{\emptyset}
\def\ts{\times}
\def\pa{\partial}                \def\os{\oplus}
\def\om{\ominus}                 \def\ev{\equiv}
\def\iint{\int\!\!\!\int}        \def\iintt{\mathop{\int\!\!\int\!\!\dots\!\!\int}\limits}
\def\el2{\ell^{\,2}}             \def\1{1\!\!1}
\def\sh{\sharp}
\def\wh{\widehat}
\def\bs{\backslash}
\def\na{\nabla}


\def\sh{\mathop{\mathrm{sh}}\nolimits}
\def\all{\mathop{\mathrm{all}}\nolimits}
\def\Area{\mathop{\mathrm{Area}}\nolimits}
\def\arg{\mathop{\mathrm{arg}}\nolimits}
\def\const{\mathop{\mathrm{const}}\nolimits}
\def\det{\mathop{\mathrm{det}}\nolimits}
\def\diag{\mathop{\mathrm{diag}}\nolimits}
\def\diam{\mathop{\mathrm{diam}}\nolimits}
\def\dim{\mathop{\mathrm{dim}}\nolimits}
\def\dist{\mathop{\mathrm{dist}}\nolimits}
\def\Im{\mathop{\mathrm{Im}}\nolimits}
\def\Iso{\mathop{\mathrm{Iso}}\nolimits}
\def\Ker{\mathop{\mathrm{Ker}}\nolimits}
\def\Lip{\mathop{\mathrm{Lip}}\nolimits}
\def\rank{\mathop{\mathrm{rank}}\limits}
\def\Ran{\mathop{\mathrm{Ran}}\nolimits}
\def\Re{\mathop{\mathrm{Re}}\nolimits}
\def\Res{\mathop{\mathrm{Res}}\nolimits}
\def\res{\mathop{\mathrm{res}}\limits}
\def\sign{\mathop{\mathrm{sign}}\nolimits}
\def\span{\mathop{\mathrm{span}}\nolimits}
\def\supp{\mathop{\mathrm{supp}}\nolimits}
\def\Tr{\mathop{\mathrm{Tr}}\nolimits}
\def\BBox{\hspace{1mm}\vrule height6pt width5.5pt depth0pt \hspace{6pt}}
\def\where{\mathop{\mathrm{where}}\nolimits}
\def\as{\mathop{\mathrm{as}}\nolimits}


\newcommand\nh[2]{\widehat{#1}\vphantom{#1}^{(#2)}}
\def\dia{\diamond}

\def\Oplus{\bigoplus\nolimits}


\newcommand{\grd}{\mathfrak{g}^r_d}
\newcommand{\gRd}{\mathfrak{g}^1_d}
\newcommand{\tC}{\widetilde{C}}
\newcommand{\tE}{\widetilde{E}}
\newcommand{\tg}{\tilde{g}}
\newcommand{\rgn}{\mathcal{R}_{g,n}}
\newcommand{\Rgn}{\overline{\mathcal{R}}_{g,n}}
\newcommand{\mg}{\mathcal{M}_g}
\newcommand{\Mg}{\overline{\mathcal{M}}_g}
\newcommand{\Mtg}{\overline{\mathcal{M}}_{\tilde{g}}}
\newcommand{\mtg}{\mathcal{M}_{\tilde{g}}}
\newcommand{\etale}{$\acute{\mbox{e}}$tale}
\newcommand{\Spec}{\mbox{Spec }}
\newcommand{\adm}{\mbox{Adm}_g(\Z_n)}
\newcommand{\pb}{\cP \Z_n}

\newcommand{\gon}{\mathop{\mathrm{gon}}}
\def\Pic{\mathop{\mathrm{Pic }}}
\def\sG{{\mathrm{G}}}

\def\qqq{\qquad}
\def\qq{\quad}
\let\ge\geqslant
\let\le\leqslant
\let\geq\geqslant
\let\leq\leqslant
\def\eq{\begin{equation}}
\def\qe{\end{equation}}
\def\bu{\bullet}

\title[{Brill-Noether theory for cyclic covers}]
{Brill-Noether theory for cyclic covers}

\date{}
\author[Irene Schwarz]{Irene Schwarz}
\address{Humboldt Universit\"at Berlin, Institut f\"ur Mathematik}


\begin{abstract}
\no We recall that the Brill-Noether Theorem gives necessary and sufficient conditions for the existence of a $\grd$.
Here we consider  a general $n$-fold,  \etale, cyclic cover $p:\tC \to C$ of a curve $C$ of genus g and investigate for which numbers $r,d$  a $\grd$ exists on $\tC$. For $r=1$ this means computing
the gonality of $\tC$. Using  degeneration to a special singular example  (containing a Castelnuovo canonical curve) and the theory of limit linear series for tree-like curves we 
show that the Pl\"ucker formula yields a necessary condition for the existence of a $\grd$ which is only slightly weaker than the sufficient condition given by the results of Laksov and Kleimann \cite{kl}, for all $n,r,d$. 
\end{abstract}

\maketitle

\section {Introduction and main results}

\setcounter{equation}{0}

We recall the Brill-Noether theorem for curves (which in this paper are always assumed to be complete reduced algebraic over $\C$). If the Brill-Noether number
\eq   \label{bnn}
\r(g,r,d):= g - (r+1)(g-d+r)
\qe 
 is non-negative,  then any smooth curve $C$ of genus $g$  carries a $\grd$, i.e. a linear
 series of dimension $r$ and degree $d$. This existence result is due to \cite{k}, \cite{kl}, see also \cite{fl}. For general curves this sufficient condition is also necessary. This is the content of the classical Brill-Noether theorem, first proved rigorously in \cite{gh}, and subsequently  in e.g.  
 \cite{eh5},  \cite{la}.

In this paper we show a Brill-Noether-like theorem  for cyclic covers
 $p:\tC\to C$, where more precisely  $p$ is a finite, flat, unramified morphism between smooth curves which has  degree $n$ and a cyclic Galois group.
We recall that, if $C$ has genus $g$,  such covers form an {\em irreducible} moduli space $\rgn$ whose birational geometry has been studied recently, see \cite{cefs}. A general point in this space is 
a cyclic cover of a general curve $C$ of genus $g$. We obtain the following non-existence result.

\begin{theorem}  \label{main}

Let $p: \tC \to C$  be a general cyclic cover of a curve of genus g. If the Brill-Noether number satisfies the inequality
\eq \label{bnest}
\r(\tg,r,d) < -r, 
\qe
where $\tg= n(g-1) +1$ is the genus of $\tC$, 
then
 $\tC$  carries no $\grd$.
\end{theorem}

To clarify the meaning of our crucial condition (\ref{bnest})  we remark that in view of (\ref{bnn}) an equivalent form is

\eq
\r(\tg - 1,r,d) < 0,
\qe
which might be understood more naturally in the context of our proof: We consider a special genus $\tg -1$ curve of compact type, obtained from normalizing a genus $\tg$ curve at one node, and then show by use of Pl\"ucker's formula that this curve already is Brill-Noether general.
We further remark that 
standard arguments immediately give the following sharpening of this theorem.

\begin{corollary}
In the situation of Theorem \ref{main}, $\tC$  carries no $\grd$, if
\eq \label{sharper}
\r(\tg,r,d) < -r+\mbox{max}\{ 0,d+1-\tg\}.  
\qe
\end{corollary}

In fact, for any divisor $D$  of degree $d$ with $h^0(D)-1=r$ on $\tC$ the divisor $K-D$ has degree $\hat{d}=2\tg-2-d$ and  $\hat{r}= h^0(K-D)-1=\tg-d+r-1$ by  Riemann-Roch.  So, if there is no $\mathfrak{g}^{\hat{r}}_{\hat{d}}$ on  $\tC$, there is no $\grd$. In particular, $\r(\tg,\hat{r},\hat{d})=\r(\tg,r,d)$. Thus, applying Theorem \ref{main} to $\hat{r}$ and $\hat{d}$, we get that the right hand side of (\ref{bnest}) could be replaced by $-\hat{r}=-r+d+1-\tg$. This proves (\ref{sharper}) and shows that the interesting part of Theorem  \ref{main} concerns small $d$.\\

Theorem \ref{main} also yields the following useful corollaries.

\begin{corollary} \label{cor}
Let $p:  \tC \to C$ be as above and assume that $r+1$ divides $\tg$. Then a general $\tC$ 
has a 
$\grd$ if and only if  
\eq
\r(\tg,r,d)\geq 0.
\qe
\end{corollary}
In fact, there exists a $d$ such that $\r(\tg,r,d)= 0$, if and only if $r+1$ divides $\tg$. For this $d$ we have $\r(\tg,r,d-k)= -k(r+1)$. So by the existence result of 
\cite{k}, \cite{kl}, there will be a $\grd$ on $\tC$, but -- by Theorem \ref{main} --  no $\mathfrak{g}^r_{d-k}$ for positive $k$. Since $\r(\tg,r,d)$
is strictly monotone in $d$, we are done.

In the case $r=1$ the  Brill-Noether theorem gives an estimate for the gonality of $\tC$, defined as
$$\gon(\tC) : =\min \{d\in \N; \mbox{ there exists a } \gRd \mbox{ on } \tC\}.$$ 
We recall from the discussion of the Brill-Noether theorem 
that in the classical language of Riemann surfaces $\gon(C)$ is the minimal number of sheets needed in a representation of $C$ as a branched cover of $\mathbb{P}^1$ which is 
$$ d = \left\lfloor \frac{g+1}{2}\right\rfloor +1 .$$

In our case of cyclic covers we have the following result:

\begin{corollary}  \label{main1}

Let $p: \tC \to C$  be a general cover as above. Then the following estimate holds

\eq \label{est}
\dfrac{\tg +1}{2} =\dfrac{n(g-1)}{2}  +1 \leq \gon(\tC) \leq \dfrac{n(g-1)}{2} + 2 = \dfrac{\tg +3}{2}.
\qe
\end{corollary}
In fact, 
for $r=1$ the smallest d such that $\tC$ has a $\gRd$ will satisfy $-1\leq\r(\tg,1,d)\leq 1$ (using Theorem \ref{main} for the lower and the existence result of \cite{k}, \cite{kl} for the upper bound). This is equivalent to the estimate (\ref{est}).

As a further consequence  of Corollary \ref{cor}, this estimate can in certain cases be sharpened.

\begin{corollary}  \label{cor1}
For $p:\tC\to C$ as above, with $n$ odd and $g$ even, we have

\eq
\gon(\tC) = \dfrac{n(g-1)+3}{2}.
\qe
\end{corollary}

The upper bound on $\gon(C)$ -- for each smooth genus g curve -- was the first rigorous result on the Brill-Noether problem (see \cite{m} with its analytic proof). The lower bound has until now only been proven by methods from algebraic geometry. Our Corollary \ref{main1} extends this result to curves which are general in the class of cyclic covers.
\\

Concerning previous results on Brill- Noether theory for cyclic covers, we recall \cite{f2}, which treats the case $r=1$ and $n=2$, i.e. gonality for cyclic double covers, and obtains a sharper result computing  $\gon(\tC)$ in each case: 
For $g$ even it is $g$, for $g$ odd it is $g+1$.
Here \cite{f2} uses Proposition 1.4.1 in \cite{f1} on 2-pointed elliptic curves.
We do not see how this argument could be mimicked even in the case $n=2, r >1$.
 This motivated our search for a different approach which finally led to this work using 
 Pl\"ucker's formula instead.
\\

Finally we emphasize that in this paper we work in characteristic 0 throughout. This allows to use the original theory of  limit linear series due to Eisenbud and Harris. We expect that our result generalizes as long as the characteristic is relatively prime to $n$. 
By \cite{dm}, this 
condition implies irreducibility of the moduli space $\rgn$  which is crucial for our proof. We remark that in arbitrary characteristic Osserman developed a (slightly) modified theory of limit linear series in \cite{o1} leading to a simple proof of the classical Brill-Noether theorem in \cite{o2}. Modifying our proof for arbitrary characteristic should either use this approach to limit linear series or generalize to arbitrary characteristic staying somewhat closer to the original approach of Eisenbud and Harris.
Formally, this is a different paper. We have not even tried.  
\\

The outline of this paper is as follows. In Section 2 we present the preliminaries needed to prove our theorem: the moduli space $\rgn$, limit linear series and Castelnuovo canonical curves.
All this (except for the use of Pl\"ucker's formula and some further results on moduli spaces of cyclic covers) is very close to the approach in 
\cite{eh5}, which simplifies the original rigorous proof of the Brill-Noether theorem given in \cite{gh}, all based on degeneration to a singular curve.
In Section 3 we prove Theorem \ref{main}.

\section {Preliminaries}
\setcounter{equation}{0}

\subsection{The moduli space $\rgn$}

In this paper we will work in the moduli space 
$\rgn$ 
of cyclic  covers as explained in the introduction.
\\
In the more traditional language of complex analysis a cyclic cover may be considered as a principal $\Z_n$-bundle, where the action is the action of the monodromy group.

By the Hurwitz formula for any cyclic cover $p: \tC \to C$ (which is unramified by definition),  the genus of $\tC$ will be $\tg=n(g-1)+1$,
if the genus of $C$ is $g$. From now on, we keep this relation between $\tg,g,n$ fixed and we use only $g,n$ as the free variables in our problem.\\

There is a bijection between  cyclic covers and {\em level $n$ curves}. See \cite{h} chapter IV exercise 2.7 for the case $n=2$, and, for general $n$,  \cite{se2} or \cite{bhpv} for a brief discussion. We need:

\begin{definition}
A \textbf{level n curve (of genus g)} or a \textbf{curve with level n structure} is a pair $(C,L)$ consisting of a smooth curve $C$ (of genus $g$) together with an $n$-torsion point $L$ of its Jacobian, i.e $L\in \Pic^0(C)$ a point of order $n$ in the group $\Pic^0(C)$.
\end{definition}

We recall that (e.g. based on the theory of stacks) one obtains a coarse moduli space $\rgn$ of level $n$ curves, see \cite{cf}. We shall need the following result:

\begin{theorem}\label{irr}
$\rgn$ is irreducible.
\end{theorem}

The proof  in the seminal paper \cite{dm} is given for the moduli space of curves with {\em full} level $n$ structure (which is smooth). But the intermediate space $\rgn$ arises as a finite quotient and thus inherits irreducibility, see \cite{b}.\\

In this paper we will also require a compactification $\Rgn$ of $\rgn$. It turns out that $\Rgn$ arises as the coarse moduli space 
of twisted level $n$ curves, defined as  a connected component of $\Mg(B\Z_n)$, which itself arises as the coarse moduli space of the category
$\overline{M}_g(B\Z_n)$ of twisted nth roots, which forms a smooth and proper Deligne-Mumford stack. 
For our notation see e.g. \cite{cefs}, using results from \cite{acv} and \cite{av}.\\

Since coarse moduli spaces are sufficient for our purpose in this paper, we shall simply take for granted the existence of $\Rgn$ as a compactification of $\rgn$, without specifying what kind of geometric objects are parametrized by this space. 
For our purpose it is sufficient to look at the closely related moduli space $\adm$ of admissible $\Z_n$-covers, which are geometrically more accessible, and describe this space in geometric terms.
We will not go into detail on the geometry of these two moduli spaces and their relation to each other. For this paper we only need to know that there is a surjective morphism $ \Rgn\to\adm$ (which in fact is the normalization map); see  Diagram \ref{diagram}.

\begin{definition}
An \textbf{admissible $\Z_n$-cover} is a morphism $p:\tC\to C$ of stable curves together with a $p$-invariant (i.e. fibre preserving) action of $\Z_n$ on $\tC$ such that:
\begin{enumerate}

\item every node of $\tC$ maps to a node of $C$
\item away from the nodes $p$ is a principal  $\Z_n$-bundle
\item for any node $x\in \tC$ there is an integer $r>0$ and a local equation $\x\e=0$ of $\tC$ at $x$, such that any element of the stabilizer of $x$ under the action of $\Z_n$ acts as $(\x,\e)\mapsto (\z\x, \z^{-1} \e)$,  where $\z$ is an rth root of unity.
\end{enumerate}
The integer r is called the \textbf{index} of x. 
 
\end{definition}

By conjugation it easily follows that every node of $\tC$ lying above the same node of $C$ has the same index $r$. Thus $r$ must divide $n$. This index can be considered as a ramification index. So an admissible  $\Z_n$-cover will be unramified, if and only if every node has index 1.\\

It is clear from the definition that for an admissible $\Z_n$-cover $p: \tC \to C$ of a smooth curve $C$, $\tC$ is also smooth and $p$ is simply a principal  $\Z_n$-bundle, which corresponds to a cyclic  cover. We shall denote the subset (in fact subscheme) of $\adm$ of such isomorphism classes of principal $\Z_n$-bundles as $\pb$.
We recall from \cite{acg}, Chapter XVI.§5, that $\pb$ actually is in bijection with 
$\rgn$, but principal $\Z_n$-bundles admit more automorphisms. These are responsible for singularities in the  moduli space $\adm$, which, however, are irrelevant in our context. \\

For the notion of families of admissible $\Z_n$-covers and an explicit construction of the coarse moduli space $\adm$ of admissible $\Z_n$-covers see \cite{acg}. This construction uses Kuranishi families of admissible $\Z_n$ covers and largely resembles the construction of the moduli space $\Mg$, avoiding the use of algebraic (or Deligne-Mumford) stacks. Thus it is both more elementary and geometrically accessible, but, if one so wishes, this moduli space can  be identified with a moduli space arising  in the algebraically more powerful theory of stacks. In the notation of \cite{cefs}, the Deligne-Mumford stack $Root_{g,n}$  of quasi-stable $n$th roots (see also \cite{j1}, \cite{j2}) is relevant. We also refer to Remark 1.6 
in \cite{cf} for the relation between stacks and admissible covers.\\

There are two forgetful maps (which are morphisms of schemes) $\adm \to\Mg$ and $\adm \to\Mtg$ that map the isomorphism class of an admissible $\Z_n$-cover $p:\tC\to C$ to the isomorphism class of  $C$ or $\tC$. Our proof only needs the map $\adm \to\Mtg$.\\

The crucial point of this very brief review is that one obtains the following commutative diagram

\begin{equation}\label{diagram}
\begin{tikzpicture} [baseline=(current  bounding  box.center)]
\matrix (m) [matrix of math nodes,row sep=2em,column sep=3em,minimum width=2em]
{\Rgn & \adm   &\Mtg \\
\rgn & \pb & \mtg \\};
     \path[->>]
    (m-1-1) edge node {} (m-1-2)
 	(m-2-1) edge node {} (m-2-2);
 	
 	\path[->]
 	(m-1-2) edge node {} (m-1-3)
 	(m-2-2) edge node {} (m-2-3);
 	
     \path[right hook->]
      (m-2-1) edge node {}(m-1-1)
      (m-2-2) edge node {} (m-1-2)
      (m-2-3) edge node {} (m-1-3);     
\end{tikzpicture}
\end{equation}

where $\twoheadrightarrow $ is a surjective morphism (as mentined above, the bottom map  actually is a bijection)  and $\hookrightarrow$ an inclusion. This is implicit in the discussion of \cite{cefs}, combining the results of Section 1.2 and
Section 1.3 (generalizing from n prime to arbitrary natural numbers $n$; for arbitrary $n$ there are additional singularities due to the factorization). This diagram is essential for our proof. On the left hand side we have irreducibility due to \cite{dm}, via stacks, and in the middle we have a geometric description which will fit with the theory of limit linear series in the following section.

\subsection{Limit linear series}

Our proof will rely largely on degenerating a family of linear series on smooth curves to a singular curve. For this we need the theory of limit linear series.

For the theory of families of linear series we refer to \cite{acg}, Chapter XXI.

We recall that (using the determinantal description of varieties of linear series) the existence of a $\grd$ is a closed condition, i.e. the set $\cG^r_d(\mg) :=\{C \in\mg\mid C \mbox{ carries a } \grd\}$ is a closed subscheme of $\mg$.

To prepare our discussion of limit linear series on a (single) tree-like curve (where the irreducible components might be nodal) we need to characterize how a single $\grd$ on a single, possibly singular, curve behaves at a smooth point.

\begin{definition} For a curve $C$, a smooth point $x\in C$ and a linear series $\ell=(L,V)\in G^r_d(C)$ we define the \textbf{vanishing sequence} of $\ell$ at x, 
$0\leq a_0(\ell, x)<a_1(\ell, x)<\cdots <a_r(\ell, x)\leq d$, as the strictly increasing order of vanishing of sections of $\ell$ in $x$, i.e. 
\eq \{a_i(\ell,x)\mid i=0,\ldots,r\} =\{ \mbox{ord}_p(\sigma )\mid \sigma \in V, \sigma \neq 0\}.\qe
The \textbf{ramification sequence} of $\ell$ at x,  $0\leq \a_0(\ell, x)\leq \a_1(\ell, x)\leq\cdots \leq \a_r(\ell, x)\leq d-r$, is defined by setting 
\eq \a_i(\ell, x):= a_i(\ell,x)-i. \qe 
The \textbf{weight} or \textbf{ramification index} of $\ell$ at x is 
\eq w^\ell(x):=\sum ^r_{i=0} \a _i(\ell, x)=\sum ^r_{i=0} a_i(\ell, x)-\frac{1}{2}r(r+1). \qe
\end{definition}

Finally, we recall the Pl\"ucker formula, which  in the following version is a crucial ingredient of our proof, see e.g. \cite{hm} Chapter 5.C Lemma 5.21 or \cite{acgh} Exercise I.C.13 for the formula and \cite{eh4} for a proof: 
\begin{theorem}    \label{PF}
If $\ell$ is any $\grd$ on a smooth and irreducible curve $C$ of genus g, then the weights will satisfy
\eq     \label{pf}
\sum_{x \in C} w^\ell(x)= (r+1) (d + r(g-1)).
\qe
\end{theorem}

Using the usual correspondence between a linear series $\ell \in  G^r_d(C)$  and a map
 $f: C \to {\mathbb P^r}$ we observe that  the ramification points of $\ell$ coincide with the ramification points of the corresponding map $f$. 
For $r=1$ the Pl\"ucker formula for $\ell$ amounts to the Riemann-Hurwitz formula for $f$.

While the definition of linear series also works for singular and reducible curves, it is useful only  in the irreducible case.
We recall from \cite{eh1} that families of $\grd$s over families of smooth curves degenerating to a reducible curve $C_0$
do not behave well with respect to limits: On $C_0$ one might obtain as limit no $\grd$ at all  or several distinct $\grd$s.
On the other hand, given a $\grd$ on the singular fibre $C_0$, it might not be possible to smooth it out, i.e. to embed it into 
a family of $\grd$s where the singular fibre is removed. 

The concept of a limit $\grd$ rectifies both problems (the second problem in full generality only for
$r=1$, for $r>1$ under some mild additional conditions). Here we shall only need  the first  (more elementary) point:  limits of linear series are limit linear series. To make the limit unique, it is very natural to consider only 1-dimensional
families of curves and associated 1-dimensional families of $\grd$s (see \cite{h}, Proposition III 9.8). However limit $\grd$s with these properties are only defined for a certain kind of nodal curves called tree-like.

We recall that a nodal curve $C$ is called {\em of compact type}, if the dual graph of C is a tree or equivalently if the Jacobian $J(C)$ is compact. A nodal curve $C$ is called {\em tree-like} if its dual graph is connected and has no cycles except, possibly, for loops.\\

The following notion of limit linear series will be central for our proof. It allows to treat singular (in particular: reducible) curves and relates them to a degeneration
through a family of smooth curves. The original theory is due to Eisenbud and Harris (see \cite{eh1}) and was originally only developed for curves of compact type. Since then it has been extended to tree-like curves (in \cite{eh2}). 

\begin{definition}
Let $C$ be a tree-like curve. A \textbf{limit linear series} (or limit $\grd$) $\ell$ of degree $d$ and dimension $r$  assigns to every irreducible component $X$ of $C$ a linear series $\ell_X \in G^r_d(X)$ called the \textbf{aspect} of $\ell$ on $X$, such that for every pair of distinct irreducible components $X,Y \subset C$ meeting in the node $x$ the aspects satisfy for all $i$
\eq 
a_i(\ell_X,x)+a_{r-i}(\ell_Y,x)\geq d.
\qe
If equality holds everywhere the limit linear series $\ell$ is called \textbf{refined}. 
\end{definition}

We recall that in \cite{eh1} our limit linear series are called {\em crude}.
The important thing about limit linear series is that they arise as limits of linear series. To make this more precise consider the following situation: Let $B$ be a regular, integral scheme of dimension 1, $b_0\in B$ a closed point and $U=B\setminus\{ b_0\}$. For example $B$ might be a smooth curve or the spectrum of a discrete valuation ring. Consider a family of stable curves $\f: \cX\to B$ such that the fibre $C_0$ over $b_0$ is tree-like and the family is smooth everywhere else, i.e. $\f|_U$ is a family of smooth curves. \\

Now let $(\cL_U, \cH_U)$ -- in the notation of \cite{acg} -- be a family of $\grd$s on $\f|_U$. 

Then for every irreducible component of $C_0$ the line bundle $\cL_U$ on $\f^{-1}(U)\subset \cX $ can be uniquely extended to a line bundle $\cL$ on $\cX$, such that $\cL|_{C_0}$ has degree $d$ on this component and degree $0$ on all other components.
Then $\cH_U$ will induce $\grd$´s on the irreducible components of $C_0$. These will turn out to form a unique limit $\grd$ on $C_0$. \\

 In \cite{eh1}, Eisenbud and Harris construct the limit $\grd$  for the case that $B$ is the spectrum of a discrete valuation ring and the special fibre is of compact type. In \cite{eh2} they extend this to a more general situation, but without giving a detailed proof. In particular, they claim:

\begin{lemma}    \label{limit}
Let $B$ be a regular, integral scheme of dimension 1, $b_0\in B$ a closed point and $U=B\setminus\{ b_0\}$. Let $\f: \cX\to B$ be a family of stable curves such that the fibre $C_0$ over $b_0$ is tree-like and the family is smooth everywhere else, i.e. $\f|_U$ is a family of smooth curves. 
If $(\cL_U, \cH_U)$ is a family of $\grd$s on $\f|_U$, then
there exists a unique crude limit $\grd$ on $C_0$ which arises as the limit of $(\cL_U, \cH_U)$.
\end{lemma}

We will need the following result, which  (if $B$ consists of more than 2 points)   is slightly stronger than the contraposition
of the statement in Lemma \ref{limit}. 
The easy special case mentioned above is the important case of $B$ being a discrete valuation ring consisting of one open and one closed point. 
We also could have reduced our proof to this special case, but we feel that the more general case is more geometrical and fits more naturally in the discussion of moduli spaces. Since we need the generalization of \cite{eh1} to tree-like curves anyway, we have formulated a general result.

\begin{proposition}\label{deg}
In the setting of Lemma \ref{limit}, if $C_0$ does not carry a crude limit $\grd$, then the general smooth fibre $C_b$ of $\f$ does not carry a $\grd$. 
\end{proposition}

The idea of a proof is as follows. Similar to $\cG^r_d(\mg)$ in $\mg$, one considers families of limit linear series
and shows that the space $$\cG^r_d(\Mg) :=\{C\in\Mg\mid C \mbox{ carries a limit-} \grd\}$$ is closed in $\Mg$. \\

Some of the most important applications of limit linear series are in using them to prove theorems on general smooth curves. 
A typical example (needed in our proof) is:
\begin{lemma}\label{EH}
A general pointed curve $(C,x)$ of genus $g$ possesses a $\grd$ with ramification sequence $0\leq \a_0(\ell, x)\leq \a_1(\ell, x)\leq\cdots \leq \a_r(\ell, x)\leq d-r$ at $x$, if and only if
\eq \label{con}
\sum_{i=0}^r \max\{\a_i+g-d+r,0\}\leq g.
\qe
\end{lemma} 

For a proof see \cite{eh2},  Proposition 1.2. The proof of sufficiency of (\ref{con}) uses the smoothing out result for limit linear series, the proof of necessity
uses Lemma \ref{limit}. Only this part is needed for our proof.

\subsection{Castelnuovo canonical curves}

We shall introduce a special class of singular curves (tree-like, but not of compact type) which will play an essential role
in our proof of Theorem \ref{main}. These Castelnuovo canonical curves go back to Castelnuovo's paper \cite{ca}, see \cite{gh} for a modern reference. 
They have topological genus 0, but arithmetical genus g. Our proof will only use the most simple case $g=1$.

\begin{definition}

An irreducible curve of arithmetical genus $g$ is called a \textbf{ Castelnuovo canonical curve} if and only if it has precisely $g$ 
simple nodes.

\end{definition}

Obviously such a curve is tree-like, since its dual graph consists of one vertex with $g$ distinct loops attached.
Its normalization is $\mathbb{P}^1$. 

Conversely, a Castelnuovo canonical curve is most easily constructed from its normalization
by selecting  $2g$ different points $x_j,y_j$ for $j=1,\cdots ,g$ and {\em identifying the points $x_j$ with $y_j$}, which is also well defined for reducible curves.
This gives  the normalization map
$$
N : \mathbb{P}^1 \to \mathbb{P}^1/ x_1 \sim y_1,\cdots, x_g \sim y_g.
$$


Since Castelnuovo canonical curves are irreducible, it makes sense to consider $\grd$s defined on them. Any such $\grd  $    $\ell$ can be pulled back to a $\grd$
on $\mathbb{P}^1 $, and one has:

\begin{lemma} \label{pullback}
The weights in smooth points coincide, i.e. if $C$ is a Castelnuovo canonical curve with normalization map $N: \mathbb{P}^1  \to C$
and $\ell$ is a $\grd$ on $C$, then 
$$ w^{N^*\ell}(x) = w^\ell(N(x))$$
for any $x \in \mathbb{P}^1 $ with $N(x)$ smooth. 
\end{lemma}

\section {Proof of Theorem \ref{main} }
\setcounter{equation}{0}

We have to show that for fixed $n,r,g,d$ and $\tg=n(g-1)+1$ 
with
\eq   \label{kp}
\r(\tg,r,d) < -r, 
\qe
 there is an open and dense set $V$ in $\rgn$ such that for any cyclic cover $p:\tC \to C$ corresponding to a point in $V$ no $\grd$ exists on $\tC$.\\

Step 1: One singular example suffices\\

First note that $\rgn$ is irreducible by Theorem \ref{irr}. Therefore every non-empty open set is already dense.
Now consider the Diagram \ref{diagram} from Section 2.1.

Since  the set
$$
\cG^r_d(\mtg)=\{C\in\mtg\mid C \mbox{ carries a } \grd\}
$$
is closed, its preimage in $\rgn$ is also closed. Therefore, its complement $V$ is open, and we are reduced to showing that it is non-empty.\\

The map $\rgn\twoheadrightarrow  \pb$ is surjective and it therefore suffices to exhibit a single principal $\Z_n$-bundle $p:\tC \to C$ in  $\pb$, such that there is no $\grd$ on $\tC$. Unfortunately, we have not succeeded in directly constructing such a smooth example. 
It is, however, possible to construct an admissible $\Z_n$-cover $p_0:\tC_0 \to C_0$ ,
such that there is no limit $\grd$ on $\tC_0$. We claim that this proves the existence of some smooth example $p:\tC \to C$, such that there is no $\grd$ on $\tC$.\\

In fact, choose a family of admissible $\Z_n$-covers over a regular one-dimensional integral base $B$, such that the special fibre over some closed point
 $b_0\in B$ corresponds to our singular example $p_0:\tC_0 \to C_0$ and the family is smooth over $B\setminus\{b_0\}$. 
This family will induce a family $\tilde{\f}:\cX\to B$ of stable curves of genus $\tg$, which satisfies the conditions of Lemma \ref{limit}. If there is no limit 
 $\grd$ on $\tC_0$ then by Proposition \ref{deg} there is some smooth fibre $\tilde{\f}^{-1}(b)=\tC$ which carries no $\grd$. 
 By the construction of $\tilde{\f}$ it is clear that this fibre comes from some principal $\Z_n$-bundle $p:\tC \to C$ . This proves our claim.\\

Step 2:  Singular examples and a preliminary estimate\\

We shall now construct  singular examples.

Fix a general 1-pointed curve $(C_0,x) \in \mathcal{M}_{g-1,1}$ and a (possibly tree-like)  1-pointed curve $(E,y) \in 
\overline{\mathcal{M}}_{1,1}$ of genus $1$, where $y$ is a smooth point of $E$. 

Fix also an admissible $\Z_n$-cover $p_E: \tE \to E$ over $E$ and set $\{y_1, \cdots, y_n\}:= p_E^{-1}(y)$. By the Hurwitz formula  $\tE$ is alsoof genus $1$.

Now choose $n$ identical copies $(C_1,x_1), \cdots (C_n,x_n)$ of $(C_0,x)$ and define the tree-like curves

\eq
C:=C_0 \cup E/ x \sim y, \qquad  \tC:= \tE \cup C_1 \cup \cdots \cup C_n / (x_1 \sim y_1, \cdots ,x_n \sim y_n).
\qe

Clearly this induces an admissible $\Z_n$-cover $p: \tC \to C$ (with index $r=1$ at each node).

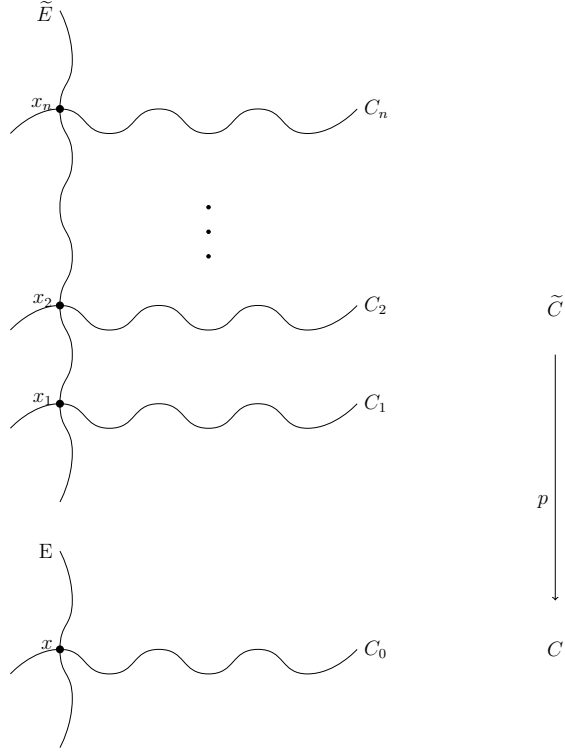
\begin{center}
\resizebox{!}{10cm}{
\begin{tikzpicture}
\draw [ xshift=4cm] plot [smooth, tension=1] coordinates { (0,0) (1/4,1) (0,2) (1/4,3)(0,4)(1/4,5)(0,6)(1/4,7)(0,8)(1/4,9)(0,10)}
node[above,left] at(0,21/10){$x_1$}
node[above,left] at(0,41/10){$x_2$}
node[above,left] at(0,81/10){$x_n$}
node[left]{$\tE$}
node at (10,4){$\tC$};

\draw[->] (14,3) to (14,-2)
node [left] at (14,0){$p$};

\draw[fill=black](4,2)circle(2pt);
\draw[fill=black](4,4)circle(2pt);
\draw[fill=black](4,8)circle(2pt);

\draw[fill=black](7,5)circle(1pt);
\draw[fill=black](7,11/2)circle(1pt);
\draw[fill=black](7,6)circle(1pt);

\draw [ xshift=4cm] plot [smooth, tension=1] coordinates {(-1,3/2)(0,2)(1,3/2)(2,2)(3,3/2)(4,2)(5,3/2)(6,2) }
node[right]{$C_1$};

\draw [ xshift=4cm] plot [smooth, tension=1] coordinates {(-1,7/2)(0,4)(1,7/2)(2,4)(3,7/2)(4,4)(5,7/2)(6,4) }
node[right]{$C_2$};

\draw [ xshift=4cm] plot [smooth, tension=1] coordinates {(-1,15/2)(0,8)(1,15/2)(2,8)(3,15/2)(4,8)(5,15/2)(6,8) }
node[right]{$C_n$};

\draw [ xshift=4cm] plot [smooth, tension=1] coordinates { (0,-5) (1/4,-4) (0,-3) (1/4,-2)(0,-1)}
node[above,left] at(0,-29/10){$x$}
node[left] {E}
node at(10,-3){$C$};

\draw[fill=black](4,-3)circle(2pt);

\draw [ xshift=4cm] plot [smooth, tension=1] coordinates {(-1,-7/2)(0,-3)(1,-7/2)(2,-3)(3,-7/2)(4,-3)(5,-7/2)(6,-3) }
node[right]{$C_0$};
\end{tikzpicture}
}
\end{center}

Assume that $\ell=\{\ell_{\tE}, \ell_{C_1}, ... , \ell_{C_n}\}$ is a limit $\grd$ on $\tC$. Denote by $$a_i: 0\leq a_{i,0}< a_{i,1} < \cdots <a_{i,r}\leq d$$ the vanishing sequence of $\ell_{\tE}$ at $x_i=y_i$,
by $$b_i: 0\leq b_{i,0}< b_{i,1}< \cdots < b_{i,r}\leq d$$ the vanishing sequence of $\ell_{C_i}$ at $x_i$
and by $$w(x_i):=w^{\ell_{\tE}}(x_i)=\sum_{j=0}^r  (a_{i,j}-j)$$ the weight of $\ell_{\tE}$ at $x_i$.
 
Then, by construction,  each $(C_i,x_i)$ is a general 1-pointed curve of genus $g-1$. In particular, it satisfies the hypothesis in Proposition \ref{EH}. Thus  the inequality (\ref{con}) holds, which implies
\eq \label{first}
\sum_{j=0}^r   b_{i,j}  \leq  (r+1)d -r(g-1) - \frac{1}{2}r(r+1).
\qe

By definition of a limit $\grd$ we have for all $i,j$

\eq \label{second}
a_{i,j}+ b_{i,r-j} \geq d.
\qe
Adding these  inequalities over all $j$ and subtracting \eqref{first} gives for all $i$

\eq \label{third}
w(x_i) \geq r(g-1).
\qe

Now assume that $E$ and $\tE$ are smooth. Then we can apply the the Pl\"ucker formula \eqref{pf} to $\tE$ to get
\footnote{If $E$ is singular, the usual Pl\"ucker formula for singular curves replaces the first equality in \ref{prel}
by $\geq $, leading to the same estimate.}
 
 \eq \label{prel}
(r+1)d = \sum_{x \in \tE} w(x) \geq \sum _{i=1}^n w(x_i) \geq  \sum _{i=1}^n r(g-1)= nr(g-1)=r(\tg -1).
 \qe 

So a $\grd$ on $\tC$ can only exist if $d\geq \dfrac{r}{r+1}(\tg -1)$. 
This is equivalent to $\r (\tg ,r,d)\geq -r(r+2)$. \\

Step 3: Improved estimate\\

It is clear from Step 1 that the estimate can be optimized by choosing $E$ and $\tE$ in a special way. 
Thus we are led to optimizing the lower bound for $d$ such that a $\grd$ with given weights $w(y_i)\geq r(g-1)$ exists on a specific elliptic curve $\tE$. 
Therefore we shall now choose $E$ and $\tE$ as Castelnuovo canonical curves.

This will allow us to lower the  left hand side of \eqref{prel} by pinching away the hole of the elliptic curve and using the Pl\"ucker formula for $g=0$.

We choose $E$ and $\tE$ as the projective line with the points $0$ and $\infty$ identified and denote its normalisation by $N:{\mathbb P^1} \to {\mathbb P^1}/0\sim \infty=E=\tE$.

The map $z\mapsto z^n$ on $\mathbb P^1$ induces an admissible $\Z_n$-cover $p_E: \tE \to E$ (with index $r=n$ at the node). As the smooth point $y$ on $E$ and its preimages $\{y_1, \cdots, y_n\}:= p_E^{-1}(y)$, we choose $y=N(1)$ and $y_i=N(\z^i)$, where $\z$ is a primitive $n$th root of unity. \\

We now have to investigate for which $g,r,d$ there exists a $\grd$ on $\tE$ with given weights $w(y_i)\geq r(g-1)$.

Note that by Lemma \ref{pullback} any $\grd \ {\tilde{\ell}}$ on $\tE$ can be pulled back to a $\grd\  \ell$ on $\mathbb{P}^1$ via the normalization $N:\mathbb{P}^1\to\tE$. Under this correspondence the weight  in any smooth point is preserved, i.e. $ w^{\tilde{\ell}}(N(x)) = w^\ell(x)$ for any point $x\in \mathbb{P}^1$ with $N(x)$ smooth.

Thus applying the Pl\"ucker formula \eqref{pf} to $\ell$ gives the following estimate:

\eq  \label{notneeded}
(r+1)(d-r) =\sum_{x \in {\mathbb P^1}} w^{\ell}(x) \geq \sum_{i=1}^n w^{\ell}(\z^i) = \sum_{i=1}^n w^{\tilde{\ell}}(y_i)\geq \sum_{i=1}^nr(g-1) = rn(g-1)=r(\tg -1).
\qe

This is equivalent to $d\geq \dfrac{r}{r+1}(\tg -1)+r$ and thus to $\r (\tg , r,d)\geq -r$, proving Theorem \ref{main}.

\section {Acknowledgements}
This paper is a restructured and substantially shortened version of my diploma thesis \cite{schw} written under the supervision of Prof. Farkas at the Humboldt University of Berlin. I want to thank my advisor G. Farkas both for suggesting this interesting subject which allowed me to work and think about a real problem in algebraical geometry already during my diploma thesis  and for contributing crucial ideas to the proof. 



\addcontentsline{toc}{section}{References}
\begin{thebibliography}{999}\setlength{\itemsep}{-\parskip} \footnotesize



%
%
%
\bibitem[ACGH]{acgh} E. Arbarello, M. Cornalba, P. A. Griffiths, J. Harris: {\em Geometry of algebraic curves}, vol. 1, Grundlehren der mathematischen Wissenschaften, Springer (1985)
%
\bibitem[ACG]{acg} E. Arbarello, M. Cornalba, P. A. Griffiths: {\em Geometry of algebraic curves},vol. 2, Grundlehren der mathematischen Wissenschaften, Springer (2011)
%
\bibitem[ACV]{acv} D. Abramovich, A. Corti, A. Vistoli: {\em Twisted bundles and admissible covers},Comm. Algebra 31 (2003), 3547--3618
%
\bibitem[AF]{f2}M. Aprodu, G. Farkas: {\em Green's conjecture for general covers}, Compact Moduli Spaces and Vector Bundles: Conference on Compact Moduli and Vector Bundles, University of Georgia, Athens, Georgia, Vol. 564,(October 21-24, 2010) . Amer. Math. Soc.(2012)
%
\bibitem[AV]{av} D. Abramovich, A. Vistoli: {\em Compactifying the space of stable maps}, J. Amer. Mat. Soc. 15 (2001), 27--75
%
\bibitem[B]{b}  M. Bernstein: {\em Moduli of curves with level structure }, Dissertation, Harvard (1999)
%
\bibitem[BHPV]{bhpv}  W. P. Barth, K. Hulek, C. A. M. Peters, A. van de Ven: {\em Compact complex surfaces}, Springer (2004)
%
%
%
%
%
\bibitem[C]{ca}   G. Castelnuovo: {\em Numero delle involutione rationali glancenti sopra una curva di dato genere}, Rend. della R. Acad. Lincei, ser. 4,5 (1889)
%
%
%
\bibitem[CEFS]{cefs}A. Chiodo, D. Eisenbud, G.Farkas, F.O. Schreyer: {\em Syzygies of torsion bundles and the geometry of the level $l$ modular variety over $\Mg$}, Invent. math. 194 (2013), 73--118
%
\bibitem[CF]{cf}A. Chiodo,  G.Farkas: {\em Singularities of the moduli space of level curves}, arXiv:1205.0201v4 (2015)
%
%
\bibitem[DM]{dm}P. Deligne, D. Mumford: {\em The irreducibility of the space of curves of given genus}, Publications Math¨¦matiques de l'IHES 36 (1969), 75--109
%
%
\bibitem[EH1]{eh4}D. Eisenbud, J. Harris: {\em Divisors on general curves and cuspidal rational curves}, Invent. math. 74 (1983), 371--418
%
\bibitem[EH2]{eh5}D. Eisenbud, J. Harris: {\em On the Brill-Noether theorem},  Algebraic geometry - open problems. Springer Berlin Heidelberg (1983) 131--137
%
\bibitem[EH3]{eh1}D. Eisenbud, J. Harris: {\em Limit linear series: basic theory}, Invent. math. 85 (1986), 337--371
%
%
\bibitem[EH4]{eh2}D. Eisenbud, J. Harris: {\em The Kodaira dimension of the moduli space of curves of genus $\leq$23}, Invent. Math. 90 (1987), 359--387
%
\bibitem[F]{f1} G. Farkas:  {\em The birational geometry of the moduli space of curves},  Academisch Proefschrift, Universiteit van Amsterdam (2000)
%
%
\bibitem[FL]{fl} W. Fulton, R. Lazarsfeld:  {\em On the connectedness of degeneracy loci and special divisors},  Acta Math. 146.1 (1981), 271-283
%
%
%
%
%
%
%
%
\bibitem[GH]{gh} P. Griffiths, J. Harris: {\em On the variety of special linear systems on a general algebraic curve},  Duke Math. J. 47(1980), 233--272
%
\bibitem[H]{h} R. Hartshorne: {\em Algebraic Geometry}, Springer (1977)
%
%
%
%
%
\bibitem[HM]{hm}J. Harris, I. Morrison: {\em Moduli of Curves}, Springer (1998)
%
%
%
%
\bibitem[J1]{j1} T.J. Jarvis: {\em Torsion-free sheaves and moduli of higher spin curves}, Compositio Mathematica 110 (1998), 291--333
%
\bibitem[J2]{j2} T.J. Jarvis: {\em Geometry of the moduli of higher spin curves}, Int. J. Math. 11 (2000), 23--47
%
\bibitem[K]{k} G. Kempf:  {\em On the geometry of a theorem of Riemann},  Ann. of Math. 98 (1973), 178--185
%
%
%
\bibitem[KL]{kl}S. Kleiman, D. Laksov:  {\em On the existence of special divisors},  Amer. J. Math. 94 (1972), 431--436
%
%
%
%
\bibitem[L]{la}  R. Lazarsfeld: {\em Brill-Noether-Petri without degeneration}, J. Diff. Geom. 23 (1986), 299 -- 307
%
\bibitem[M]{m} Th. Meis:  {\em Die minimale Bl\"atterzahl der Konkretisierungen einer kompakten Riemannschen Fl\"ache}, Schriftenreihe d. Math. Inst. d. Univ. M\"unster, H. 16 (1960)
%
%
%
\bibitem[O1]{o1} B. Ossermann: {\em A limit linear series moduli scheme}, Annales de l'Institut Fourier, 56 (2006), 1165 --1205
%
\bibitem[O2]{o2} B. Ossermann: {\em A simple characteristic-free proof of the Brill-Noether theorem}, Bull. Braz. Math. Soc., New Series 45(4) (2014), 807 --818
%
%
%
%
%


%
%
\bibitem[S]{se2} J.P. Serre: {\em Sur la topologie des varietes algebriques en caracteristique p}, Symposium de topologie algebrique, Mexico (1956), 24--53

\bibitem[Sch]{schw} I. Schwarz: {\em Brill-Noether theory for cyclic covers}, diploma thesis, Berlin (2016)
%
%
%
%
%
\end{thebibliography}
\end{document}